\documentclass[12pt, a4paper, twoside]{article}
\usepackage{amssymb}
\usepackage[leqno]{amsmath}
\usepackage{latexsym}
\usepackage{enumerate}
\usepackage{amsmath, amssymb, amscd}
\Roman{equation}

\DeclareMathOperator{\Pic}{Pic}
\DeclareMathOperator{\Supp}{Supp}

\begin{document}

\title{\textbf{\Large{New examples of Weierstrass semigroups associated with a double covering of a curve on a Hirzebruch surface of degree one}}}

\author{Kenta Watanabe \thanks{Nihon University, College of Science and Technology,   7-24-1 Narashinodai Funabashi city Chiba 274-8501 Japan , {\it E-mail address:watanabe.kenta@nihon-u.ac.jp}, Telephone numbers: 090-9777-1974} }

\date{}

\maketitle 


\begin{abstract}

\noindent Let $\varphi:\Sigma_1\longrightarrow \mathbb{P}^2$ be a blow up at a point on $\mathbb{P}^2$. Let $C$ be the proper transform of a smooth plane curve of degree $d\geq 4$ by $\varphi$, and let $P$ be a point on $C$.  Let $\pi:\tilde{C}\longrightarrow C$  be a double covering branched along the reduced divisor on $C$ obtained as the intersection of $C$ and a reduced divisor in $|-2K_{\Sigma_1}|$ containing $P$. In this paper, we investigate the Weierstrass semigroup $H(\tilde{P})$ at the ramification point $\tilde{P}$ of $\pi$ over $P$, in the case where the intersection multiplicity at $\varphi(P)$ of $\varphi(C)$ and  the tangent line at $\varphi(P)$ of $\varphi(C)$ is $d-1$.

\end{abstract}

\section{Introduction}

\noindent We work over the complex number field $\mathbb{C}$. Let $C$ be a smooth projective curve, and let $P$ be a point on $C$. Then we call a natural number $n$ satisfying 
$$h^0(\mathcal{O}_C((n-1)P))=h^0(\mathcal{O}_C(nP))$$
a {\it{gap}} at $P$ on $C$. Let $\mathbb{N}_0$ be the additive monoid consisting of non-negative integers. If we let $G(P)$ be the set of gaps at $P$ on $C$, the set $\mathbb{N}_0\backslash G(P)$ forms an additive monoid. We call it the {\it{Weierstrass semigroup}} at $P$ on $C$, and denote it by $H(P)$. If the genus of $C$ is $g$, then the cardinality of $G(P)$ coincides with $g$ and it is called the genus of $H(P)$.

Assume that $C$ is a double covering of $\mathbb{P}^1$, that is, $C$ is hyperelliptic. Then the ramification points of it are Weierstrass points. Hence, in general,  it is natural and interesting to consider the problem of whether the ramification points of a given double covering of a curve are Weierstrass.

Let $\pi:\tilde{C}\longrightarrow C$ be a double covering, let $P$ be a branch point of it, and let $\tilde{P}$ be the ramification point of $\pi$ over $P$. Then we consider the Weierstrass semigroup $H(\tilde{P})$ called the {\it{double covering type}}. If the genus of $\tilde{C}$ is $\tilde{g}$ and $C=\mathbb{P}^1$, $H(\tilde{P})$ is generated by 2 and $2\tilde{g}+1$. Komeda [2] has determined the all possible Weierstrass semigroups $H(\tilde{P})$, in the case where the genus $g$ of $C$ is 3 and the genus of  $\tilde{C}$ is $\tilde{g}\geq9$. On the other hand, in the case where $C$ is a smooth plane curve of degree $d\geq4$  and the intersection multiplicity at $P$ of $C$ and the tangent line at $P$ of $C$ is $d$ or $d-1$, we have determined $H(\tilde{P})$, under the condition that the branch locus of $\pi$ is the reduced divisor on $C$ of degree $6d$ obtained as the intersection of $C$ and a reduced divisor on $\mathbb{P}^2$ of degree 6 ([4] and [5]). If $\tilde{C}$ is smooth, that is, the branch divisor of $\pi$ is reduced, the computations of  the Weierstrass semigroups at the ramification points are interesting from the point of view of the classification of Weierstrass points on a curve. However, in this paper, we will focus on the Weierstrass semigroups at the points on a curve obtained by the normalization of $\tilde{C}$, in the case where the branch divisor of $\pi$ is not reduced.

Let $\varphi:\Sigma_1\longrightarrow \mathbb{P}^2$ be a blow up at a point on $\mathbb{P}^2$. Then we call the surface $\Sigma_1$ a Hirzebruch surface. Let $E$ be the exceptional divisor of $\varphi$, and we denote by $L$ the total transform of a line on $\mathbb{P}^2$ by $\varphi$.  Then we note that the Picard lattice $\Pic(\Sigma_1)$ of $\Sigma_1$ is generated by the classes of $E$ and $L$. In the previous work, we have the following result.

\newtheorem{thm}{Theorem}[section]

\begin{thm} {\rm{([6, Theorem 1.2]).}} Let $C$ be a smooth curve on $\Sigma_1$  which is linearly equivalent to the divisor $dL-E$ with $d\geq4$. Let $\pi:\tilde{C}\longrightarrow C$ be a double covering of $C$ with a branch point $P$, $\pi^{-1}(P)=\{\tilde{P}\}$, and assume that there exists an effective divisor $T_P\in |L|$ such that $T_P|_C=dP$. Moreover, let $M_d$ be the following condition.

\smallskip

\smallskip

\noindent $M_d$ : There exists a double covering $\tilde{\pi}:X\longrightarrow \Sigma_1$ branched along a reduced divisor belonging to $|-2K_{\Sigma_1}|$ such that $\tilde{C}\subset X$ and $\tilde{\pi}|_{\tilde{C}}=\pi$.

\smallskip

\smallskip

\noindent Then, we have the following results.

\smallskip

\smallskip

\noindent {\rm{(a)}} If $P\in E$, then the following conditions are equivalent.

\smallskip

\smallskip

{\rm{(i)}} $M_d$ holds.

\smallskip

{\rm{(ii)}} $H(\tilde{P})=2H(P)+(6d-3)\mathbb{N}_0$.

\smallskip

\smallskip

\noindent {\rm{(b)}} If $P\notin E$ and $M_d$ holds, then $H(\tilde{P})=2H(P)+(6d-1)\mathbb{N}_0+(2d^2+1)\mathbb{N}_0.$\end{thm}

\noindent In Theorem 1.1, if the double covering $\pi$ satisfies the condition $M_d$,  $\tilde{C}$ is the normalization of the double covering of $\varphi(C)$ branched along the divisor obtained as the restriction to $\varphi(C)$ of a reduced divisor on $\mathbb{P}^2$ of degree 6 which has a singularity at $\varphi(E)$. However,  conversely, it is difficult to consider the problem of whether $\pi$ satisfies the condition $M_d$, even if the Weierstrass semigroup $H(\tilde{P})$ at the ramification point $\tilde{P}$ of $\pi$ over the branch point $P$ is given. In this paper, we will investigate the Weierstrass semigroup $H(\tilde{P})$, in the case where $C$ is the proper transform of a smooth plane curve of degree $d\geq4$ by $\varphi$, and the intersection multiplicity at $\varphi(P)$ of $\varphi(C)$ and the tangent line at $\varphi(P)$ of $\varphi(C)$ is $d-1$. Moreover, we consider a necessary and sufficient condition for a double covering $\pi$ of a smooth curve on a Hirzebruch surface to satisfy the condition $M_d$. Our main theorem is as follows.

\begin{thm} Let $L$, $E$, and $C$ be as in Theorem 1.1. Let $\pi:\tilde{C}\longrightarrow C$ be a double covering which has a branch point $P$ on $C$, and let  $\pi^{-1}(P)=\{\tilde{P}\}$. Assume that there exist $T_P\in |L|$, point $Q\in C$ with $Q\neq P$, and $T_Q\in |L|$ such that ${T_P}|_C=(d-1)P+Q$ and  ${T_Q}|_C=dQ$, and let $M_d$ be the following condition.

\smallskip

\smallskip

\noindent $M_d$ : There exists a double covering $\tilde{\pi}:X\longrightarrow \Sigma_1$ branched along a reduced divisor belonging to $|-2K_{\Sigma_1}|$ such that $\tilde{C}\subset X$ and $\tilde{\pi}|_{\tilde{C}}=\pi$. 

\smallskip

\smallskip

\noindent Then we get the following results.

\smallskip

\smallskip

\noindent {\rm{(a)}} If  $Q\in E$, then the following conditions are equivalent.

\smallskip

{\rm{(i)}} $\pi$ satisfies $M_d$.

\smallskip

{\rm{(ii)}} $H(\tilde{P})=2H(P)+(8d-9)\mathbb{N}_0+\cdots+(8d-9+2r(d-2))\mathbb{N}_0$

$+\cdots+(8d-9+2(d-3)(d-2))\mathbb{N}_0$.

\smallskip

\smallskip

\noindent {\rm{(b)}} If $P\in E$, then we get the following assertion.

\smallskip

\smallskip

{\rm{(i)}} If $\pi$ satisfies $M_d$, then $H(\tilde{P})= 2H(P)+(8d-11)\mathbb{N}_0$

$+\cdots+(8d-11+2r(d-2))\mathbb{N}_0+\cdots+(8d-11+2(d-4)(d-2))\mathbb{N}_0$.

\smallskip

{\rm{(ii)}} Assume that $Q$ is a branch point of $\pi$, and let $\tilde{Q}$ be the ramification point over $Q$. If $2d^2-3\notin H(\tilde{P})$ and  $6d-1\in H(\tilde{Q})$, then $\pi$ satisfies $M_d$.

\smallskip

\smallskip

\noindent {\rm{(c)}} If $P\notin E$, $Q\notin E$, and $\pi$ satisfies $M_d$, then $H(\tilde{P})=2H(P)+(8d-9)\mathbb{N}_0$

$+\cdots+(8d-9+2r(d-2))\mathbb{N}_0+\cdots+(8d-9+2(d-4)(d-2))\mathbb{N}_0+(2d^2-1)\mathbb{N}_0$.

 \end{thm}

\noindent In Theorem 1.2 (b), if $H(\tilde{P})$ is the Weierstrass semigroup as in (i), then $2d^2-3\notin H(\tilde{P})$. This means that the assertion of (ii) gives the converse assertion of (i) and Theorem 1.1 (b).

$\;$

\noindent {\bf{Notations and Conventions}}. A curve and a surface are smooth and projective. For a curve or a surface $Y$, we denote a canonical divisor of $Y$ by $K_Y$. If $C$ is a curve on a surface $Y$, by the adjunction formula, $K_C=(K_Y+C)|_C$. For a curve $C$, we denote the genus of $C$ by $g(C)$. For a divisor $D$ on a curve or a surface, we denote by $|D|$ the linear system defined by it.  If two divisors $D_1$ and $D_2$ belong to the same linear system, we will write $D_1\sim D_2$. For two divisors $D_1$ and $D_2$ on a surface $Y$ and $R\in \Supp D_1\cap\Supp D_2$, we denote by $I_R(D_1\cap D_2)$ the intersection multiplicity at $R$ of  $D_1$ and $D_2$. We call the minimum degree of pencils on a curve $C$ the {\it{gonality}} of $C$. 

A submonoid $H\subset \mathbb{N}_0$ is called a {\it{numerical semigroup}} if the set $\mathbb{N}_0\backslash H$ is a finite set. The {\it{genus}} of a numerical semigroup $H$ is defined by the cardinality of $\mathbb{N}_0\backslash H$, and it is denoted by $g(H)$. For a numerical semigroup $\tilde{H}$, if we let 
$$d_2(\tilde{H}):=\{\frac{h}{2}\;|\;h\text{ is even and } h\in \tilde{H}\},$$
\noindent it is a numerical semigroup. 

\section{Proof of Theorem 1.2}

In this section, we will give a proof of Theorem 1.2 and some examples of it. Assume that the notations are as in Theorem 1.2. Then $C$ is the proper transform of a smooth curve of degree $d\geq4$ by the blow up $\varphi:\Sigma_1\longrightarrow \mathbb{P}^2$ at a point on $\mathbb{P}^2$. Moreover, it is well known that $d_2(H(\tilde{P}))=H(P)$ (cf. [3]). Hence, the results on the computation of the Weierstrass semigroup at a Weierstrass point on a smooth plane curve are often used to compute the Weierstrass semigroup $H(\tilde{P})$. First of all, we recall the following useful result.

\newtheorem{lem}{Lemma}[section]

\begin{lem} {\rm{(cf. [1] TABLE 3 and TABLE 4)}}. Let $C$ be a plane curve of degree $d\geq 4$, and let $P\in C$. Let $T_P$ be the tangent line of $C$ at $P$. Then, we have the following results.

\smallskip

\smallskip

\noindent {\rm{(i)}} If $I_P(C\cap T_P)=d$, then $H(P)=d\mathbb{N}_0+(d-1)\mathbb{N}_0$.

\noindent {\rm{(ii)}} If $I_P(C\cap T_P)=d-1$, then $H(P)=(d-1)\mathbb{N}_0+(2d-3)\mathbb{N}_0$

$+\cdots+(d-1+r(d-2))\mathbb{N}_0+\cdots+(d-1+(d-2)^2)\mathbb{N}_0$. \end{lem}

\smallskip

\smallskip

\noindent {\it{Proof of Theorem 1.2}}. Let  $\tilde{\pi}:X\longrightarrow \Sigma_1$ be a double covering of $\Sigma_1$ as in the condition $M_d$. Since $K_{\Sigma_1}\sim -3L+E$ and $\tilde{C}$ is smooth, $C$ intersects transversely the branch locus of $\tilde{\pi}$ at distinct $6d-2$ smooth points of it. Let $\eta:\tilde{X}\longrightarrow X $ be a minimal resolution of $X$. Since $\tilde{C}$ is smooth, it does not contain any singular point of $X$. Hence, $\eta^{-1}(\tilde{C})=\tilde{C}$. Since the exceptional divisor of $\eta$ does not intersect $\eta^{-1}(\tilde{C})$, we have $K_{\tilde{C}}\sim\tilde{C}|_{\tilde{C}}$. From now on let $\tilde{D}:=(\tilde{\pi}\circ\eta)^{\ast}D$, for a divisor $D$ on $\Sigma_1$.

\smallskip

\smallskip

(a) (i) $\Longrightarrow$ (ii) We classify divisors $D$ on $\Sigma_1$ which are linearly equivalent to $C$ and $D|_C$ are effective to investigate the set of gaps $G(\tilde{P})$ at $\tilde{P}$. Let $L_P\in |L|$ be a divisor such that $P\in L_P$ and $L_P\neq T_P$. For $0\leq s\leq d$ and $0\leq t\leq\min\{d-2,d-s\}$, let $D_1=sT_P+tL_P+(d-s-t)T_Q-E$. Moreover, let $D_2=(d+1)T_P-T_Q-E$. Since  $I_{\tilde{P}}(\tilde{D_1}\cap\tilde{C})=2s(d-1)+2t$, we have $2s(d-1)+2t+1\in\mathbb{N}_0\backslash H(\tilde{P})$. Since $K_{\tilde{C}}=\tilde{D_2}|_{\tilde{C}}=(2d^2-2)\tilde{P}$, $2d^2-1\in \mathbb{N}_0\backslash H(\tilde{P})$. Therefore, we have
$$\mathbb{N}_0\backslash H(\tilde{P})=\{2s(d-1)+2t+1\;|\;0\leq s\leq d,\; 0\leq t\leq\min\{d-2,d-s\}\}\cup\{2d^2-1\}.$$
Indeed, the cardinality of the set of the right hand side is $\displaystyle\frac{d^2+3d-2}{2}$ and it coincides with $g(\tilde{C})-g(C)$.
If $2s(d-1)+2(d-s)+3\leq2(s+1)(d-1)-1$, then $s\geq3$. Hence, the minimum odd number of $H(\tilde{P})$ is $8d-9$. We have $2H(P)+(8d-9)\mathbb{N}_0\subset H(\tilde{P})$. Let $n$ be an odd number satisfying $n\in H(\tilde{P})\backslash (2H(P)+(8d-9)\mathbb{N}_0)$. Then we have $n\leq 8d-9+2(2g(C)-1)=2d^2+2d-7$.

Assume that $2d^2+1\leq n\leq 2d^2+2d-7$. Then there exists an odd number $k$ such that $3\leq k\leq 2d-5$ and $n=2(d+2)(d-1)-k$. We have 
$$n-(8d-9)=(2d-3-k)(d-1)+(k-1)(d-2)\in 2H(P).$$
This is a contradiction. 

Assume that $2d^2-2d+3\leq n\leq 2d^2-3$. Then there exists an odd number $k$ such that $1\leq k\leq 2d-5$ and $n=2(d^2-1)-k$. Since 
$$n-(8d-9)=(2d-5-k)(d-1)+(k-1)(d-2)\notin 2H(P),$$
we have $k=2d-5$. Hence, we have $n=8d-9+2(d-3)(d-2)=2d^2-2d+3$. 

Assume that there exists an integer $s$ satisfying 
$$3\leq s\leq d-1\text{ and }2s(d-1)+2(d-s)+3\leq n\leq 2(s+1)(d-1)-1.$$
Then there exists an odd number $k$ with $1\leq k\leq 2s-5$ and $n=2(s+1)(d-1)-k$. Since $n-(8d-9)=(2s-5-k)(d-1)+(k-1)(d-2)\notin 2H(P)$, we have $k=2s-5$. This means that there exists an integer $r$ with $0\leq r\leq d-4$ and $n=8d-9+2r(d-2)$.
Therefore, we have the assertion. 

\smallskip

\smallskip

(ii) $\Longrightarrow$ (i) Let $B$ be the branch divisor of $\pi$. Then by the Hurwitz formula, there exists an effective divisor $D$ of degree $3d-1$ on $C$ with $B\sim 2D$. Let $D^{'}=-K_{\Sigma_1}|_C-D$. Since $T_P|_C\sim T_Q|_C$, we have $(d-1)P\sim (d-1)Q$. Since $C$ is linearly equivalent to $dT_P-E$ as a divisor on $\Sigma_1$, if we let $F=K_C+D^{'}+P$, we have $F\sim d^2P-D$. Since $\deg(D^{'})=0$, we show that $h^0(\mathcal{O}_C(-D^{'}))>0$ to prove that $B\in|-2K_{\Sigma_1}|_C|$. Assume that $h^0(\mathcal{O}_C(-D^{'}))=0$. Since $K_C-F+P=-D^{'}$, we have $h^0(\mathcal{O}_C(K_C-F+P))=h^0(\mathcal{O}_C(K_C-F))=0$. This means that $h^0(\mathcal{O}_C(d^2P-D))=h^0(\mathcal{O}_C((d^2-1)P-D))+1$. On the other hand, since $d^2\in H(P)$, we have $h^0(\mathcal{O}_C(d^2P))=h^0(\mathcal{O}_C((d^2-1)P))+1$. Since 
$$h^0(\mathcal{O}_{\tilde{C}}(2d^2\tilde{P}))=h^0(\mathcal{O}_C(d^2P))+h^0(\mathcal{O}_C(d^2P-D))$$
and $h^0(\mathcal{O}_{\tilde{C}}((2d^2-2)\tilde{P}))=h^0(\mathcal{O}_C((d^2-1)P))+h^0(\mathcal{O}_C((d^2-1)P-D))$, we have $h^0(\mathcal{O}_{\tilde{C}}(2d^2\tilde{P}))=h^0(\mathcal{O}_{\tilde{C}}((2d^2-2)\tilde{P}))+2$. This contradicts the fact that $2d^2-1\notin H(\tilde{P})$. Hence, $\pi$ satisfies $M_d$.

$\;$

(b) (i) By the same way as above, we determine the set of gaps $G(\tilde{P})$ at $\tilde{P}$. Let $L_P$ be as above. For $0\leq s\leq d$ and $0\leq t\leq\min\{d-2,d-s\}$ with $(s,t)\neq(0,0)$, we set $D_1=sT_P+tL_P+(d-s-t)T_Q-E$. Moreover, let $D_2=(d+1)T_P-T_Q-E$ and $D_3=L_P+dT_P-T_Q-E$. Then $\tilde{D_i}\sim \tilde{C}$ and $D_i|_C$ is effective. First of all, we have $I_{\tilde{P}}(\tilde{D_1}\cap\tilde{C})=2s(d-1)+2t-2$. Hence, $2s(d-1)+2t-1\in \mathbb{N}_0\backslash H(\tilde{P})$. Since $I_{\tilde{P}}(\tilde{D_2}\cap\tilde{C})=2d^2-4$, we have $2d^2-3\in \mathbb{N}_0\backslash H(\tilde{P})$. Since $I_{\tilde{P}}(\tilde{D_3}\cap\tilde{C})=2d^2-2d$, we have $2d^2-2d+1\in \mathbb{N}_0\backslash H(\tilde{P})$. By the computation of the cardinality of $\mathbb{N}_0\backslash H(\tilde{P})$, we have

\smallskip

\smallskip

\noindent $\mathbb{N}_0\backslash H(\tilde{P})=\{2s(d-1)+2t-1\;|\;0\leq s\leq d,\; 0\leq t\leq\min\{d-2,d-s\}, \text{ and }(s,t)\neq(0,0)\}\cup\{2d^2-2d+1,2d^2-3\}.$

\smallskip

\smallskip

If $2s(d-1)+2(d-s)+1\leq2(s+1)(d-1)-3$, then $s\geq3$. Hence, the minimum odd number of $H(\tilde{P})$ is $8d-11$. Therefore, we have $2H(P)+(8d-11)\mathbb{N}_0\subset H(\tilde{P})$. Let $n$ be an odd number such that $n\in H(\tilde{P})\backslash (2H(P)+(8d-11)\mathbb{N}_0)$. Then we have $n\leq 8d-11+2(2g(C)-1)=2d^2+2d-9$. 
We note that there exists an integer $s$ satisfying  $3\leq s\leq d-1$, and $2s(d-1)+2(d-s)+1\leq n\leq 2(s+1)(d-1)-3$. 
In fact, if we assume that $2d^2-2d+3\leq n\leq 2d^2-5$, there exists an odd number $k$ such that $3\leq k\leq 2d-5$ and $n=2(d^2-1)-k$. Hence, we have
$$n-(8d-11)=(2d-3-k)(d-1)+(k-3)(d-2)\in 2H(P).$$
This is a contradiction. If we assume that $2d^2-1\leq n\leq 2d^2+2d-9$, then there exists an odd number $k$ such that $5\leq k\leq 2d-3$ and $n=2(d+2)(d-1)-k$.  Hence, we have
$$n-(8d-11)=(2d-1-k)(d-1)+(k-3)(d-2)\in 2H(P).$$
This is a contradiction. Therefore, there exists an odd number $k$ such that
$$3\leq k\leq 2s-3\text{ and }n=2(s+1)(d-1)-k.$$
Since $n-(8d-11)=(2s-3-k)(d-1)+(k-3)(d-2)\notin 2H(P)$, we have $k=2s-3$. Therefore, there exists an integer $r$ such that $0\leq r\leq d-4$ and $n=8d-11+2r(d-2)$. Hence, we have the assertion.

\smallskip

\smallskip

(ii) Let $D$ be a divisor with $B\sim 2D$ for the branch divisor $B$ of $\pi$. Let $D^{'}=-K_{\Sigma_1}|_C-D-Q$. Since $E\cap C=\{P\}$, by the same reason as above, if we let $F=K_C+D^{'}+P$, we have $F\sim(d^2-1)P-D$. Then we have $h^0(\mathcal{O}_C(-D^{'}))>0$. Assume that $h^0(\mathcal{O}_C(-D^{'}))=0$. Then $h^0(\mathcal{O}_C(K_C-F))=0$. This means that $h^0(\mathcal{O}_C((d^2-1)P-D))=h^0(\mathcal{O}_C((d^2-2)P-D))+1$. On the other hand, since $d^2-1\in H(P)$, we have $h^0(\mathcal{O}_C((d^2-1)P))=h^0(\mathcal{O}_C((d^2-2)P))+1$. Since $h^0(\mathcal{O}_{\tilde{C}}((2d^2-2)\tilde{P}))=h^0(\mathcal{O}_C((d^2-1)P))+h^0(\mathcal{O}_C((d^2-1)P-D))$ and $h^0(\mathcal{O}_{\tilde{C}}((2d^2-4)\tilde{P}))=h^0(\mathcal{O}_C((d^2-2)P))+h^0(\mathcal{O}_C((d^2-2)P-D))$, we have $h^0(\mathcal{O}_{\tilde{C}}((2d^2-2)\tilde{P}))=h^0(\mathcal{O}_{\tilde{C}}((2d^2-4)\tilde{P}))+2$. However, this contradicts the assumption that $2d^2-3\notin H(\tilde{P})$. Since $\deg(-D^{'})=1$, there exists a point $Q^{'}$ on $C$ belonging to $|-D^{'}|$.  Then we have 
$$D\sim -K_{\Sigma_1}|_C-Q+Q^{'}.\leqno (2.1)$$

On the other hand, since $T_Q|_C=dQ$, by Lemma 2.1 (i), we have $3d\in H(Q)$. Hence, we have $6d\in H(\tilde{Q})$. Moreover, by the assumption that $6d-1\in H(\tilde{Q})$, we have $h^0(\mathcal{O}_{\tilde{C}}(6d\tilde{Q}))=h^0(\mathcal{O}_{\tilde{C}}((6d-2)\tilde{Q}))+2$. Here, we set $D^{''}\sim 3dQ-D$. Since $h^0(\mathcal{O}_{\tilde{C}}(6d\tilde{Q}))=h^0(\mathcal{O}_{C}(3dQ))+h^0(\mathcal{O}_{C}(D^{''}))$ and
$$h^0(\mathcal{O}_{\tilde{C}}((6d-2)\tilde{Q}))=h^0(\mathcal{O}_{C}((3d-1)Q))+h^0(\mathcal{O}_{C}(D^{''}-Q)),$$ we have
$$h^0(\mathcal{O}_{C}(D^{''}))=h^0(\mathcal{O}_{C}(D^{''}-Q))+1>0.\leqno (2.2)$$
Since $\deg(D^{''})=1$, there exists a point $Q^{''}$ on $C$ belonging to $|D^{''}|$. Since, by the linear equivalence (2.1), 
$$3dQ-Q^{''}\sim D \sim (3T_Q-E)|_C-Q+Q^{'}=3dQ-P-Q+Q^{'},$$
we have $P+Q\sim Q^{'}+Q^{''}$. Since $C$ is isomorphic to the smooth plane curve $\varphi(C)$ via $\varphi$, the gonality of $C$ is $d-1$ and the divisor $(T_P-E)|_C=(d-2)P+Q$ gives a gonality pencil on $C$. Hence, $\dim |P+Q|=0$. This means that $P+Q=Q^{'}+Q^{''}$. By the equality (2.2), we have $Q^{''}\neq Q$, and hence, we have $Q^{'}=Q$ and $Q^{''}=P$. Therefore, we have $D\sim -K_{\Sigma_1}|_C$ and hence, $\pi$ satisfies $M_d$.

$\;$

(c) Let $C\cap E=\{R\}$ and let $L_{Q,R}, L_{P,R}\in |L|$ be divisors such that $Q,R\in L_{Q,R}$ and $P,R\in L_{P,R}$. For integers $s$ and $t$ satisfying $0\leq s\leq d$, $0\leq t\leq \min\{d-2,d-s\}$, and $(s,t)\neq(d,0)$ let 
$$D_1=sT_P+tL_{P,R}+(d-s-t)L_{Q,R}-E.$$
Moreover, let $D_2=L_{P,R}+dT_P-T_Q-E$ and let $D_3=L_{Q,R}+dT_P-T_Q-E$. Since $I_{\tilde{P}}(\tilde{D_1}\cap\tilde{C})=2s(d-1)+2t$, we have $2s(d-1)+2t+1\in\mathbb{N}_0\backslash H(\tilde{P})$. Since $I_{\tilde{P}}(\tilde{D_2}\cap\tilde{C})=2d(d-1)+2$, we have $2d^2-2d+3\in\mathbb{N}_0\backslash H(\tilde{P})$. Since $I_{\tilde{P}}(\tilde{D_3}\cap\tilde{C})= 2d(d-1)$, we have $2d(d-1)+1\in\mathbb{N}_0\backslash H(\tilde{P})$. Therefore, by the same reason as above, we have
$$\mathbb{N}_0\backslash H(\tilde{P})=\{2s(d-1)+2t+1\;|\;0\leq s\leq d,\; 0\leq t\leq\min\{d-2,d-s\}\}\cup\{2d^2-2d+3\}.$$
If $2s(d-1)+2(d-s)+3\leq2(s+1)(d-1)-1$, then $s\geq3$. Hence, the minimum odd number of $H(\tilde{P})$ is $8d-9$. We have $2H(P)+(8d-9)\mathbb{N}_0\subset H(\tilde{P})$. Let $n$ be an odd number satisfying $n\in H(\tilde{P})\backslash (2H(P)+(8d-9)\mathbb{N}_0)$. Then we have $n\leq 8d-9+2(2g(C)-1)=2d^2+2d-7$. By the same reason as in the proof of (a), the case that $2d^2+1\leq n\leq 2d^2+2d-7$ does not occur. 

Assume that $2d^2-2d+5\leq n\leq 2d^2-1$. Then there exists an odd number $k$ such that $-1\leq k\leq 2d-7$ and $n=2(d^2-1)-k$. Since 
$$n-(8d-9)=(2d-5-k)(d-1)+(k-1)(d-2)\notin 2H(P),$$
we have $k=-1$. Hence, we have $n=2d^2-1$. 

Assume that there exists an integer $s$ satisfying 
$$3\leq s\leq d-1\text{ and }2s(d-1)+2(d-s)+3\leq n\leq 2(s+1)(d-1)-1.$$
By the same reason as in the proof of (a), there exists an integer $r$ such that $0\leq r\leq d-4$ and $n=8d-9+2r(d-2)$. Therefore, we have the assertion. $\hfill\square$

\smallskip

\smallskip

\noindent We can construct an example for each case as in Theorem 1.2.

$\;$

\noindent {\bf{Example 2.1}}. Let $B_0\subset\mathbb{P}^2$ be a reduced and irreducible divisor of degree 6 defined by an equation $f(x,y,z)=0$. Assume that $B_0$ has a double point $R$, and $B_0$ does not have other singularity. Let $X_0$ be the hypersurface in $\mathbb{P}(1,1,1,3)$ defined by the equation $w^2=f(x,y,z)$. Let $X$ be the blow up at the point on $X_0$ corresponding to $R$. If $\varphi:\Sigma_1\longrightarrow \mathbb{P}^2$ is the blow up at $R$ and $B$ is the proper transform of $B_0$ by $\varphi$, $X$ is obtained as the double covering $\tilde{\pi}:X\longrightarrow\Sigma_1$ of $\Sigma_1$ branched along $B$.

Let $C_0\subset\mathbb{P}^2$ be the plane curve defined by the equation $yz^{d-1}-x^{d-1}z+y^d=0$. Assume that $R\in C_0$ and $C_0$ intersects transversely $B_0$ at other smooth points of $B_0$. Let $C$ be the proper transform of $C_0$ by $\varphi$. Moreover, we set $\varphi^{-1}(R)=E$. We note that if we let $L\in |\varphi^{\ast}\mathcal{O}_{\mathbb{P}^2}(1)|$ be a divisor, then $C\in |dL-E|$. Let $\tilde{C}:=\tilde{\pi}^{-1}(C)$. Then $\tilde{\pi}$ induces the double covering $\pi:\tilde{C}\longrightarrow C$ branched along the divisor $B\cap C$ on $C$. Let $\varphi^{-1}(0:0:1)\cap C=\{P\}$. Assume that $P$ is a branch point of $\pi$ and let $\pi^{-1}(P)=\{\tilde{P}\}$. 

\smallskip

\smallskip

(a) We consider the case where $R=(1:0:0)$.  Let $E\cap C=\{Q\}$. Since $T_R|_{C_0}=dR$ and $T_{\varphi(P)}|_{C_0}=(d-1)\varphi(P)+R$, if we let $T_Q=\varphi^{-1}(T_R)$ and $T_P=\varphi^{-1}(T_{\varphi(P)})$, we have $T_Q|_C=dQ$ and $T_P|_C=(d-1)P+Q$.  Hence, $H(\tilde{P})$ coincides with the Weierstrass semigroup as in Theorem 1.2 (a) (ii).

\smallskip

(b) We consider the case where $R=(0:0:1)$. Then $P\in E$. We note that if $I_R(C_0\cap B_0)=3$, this case can occur. Let $\varphi^{-1}(1:0:0)=\{Q\}$. Since $T_R|_{C_0}=(d-1)R+\varphi(Q)$ and $T_{\varphi(Q)}|_{C_0}=d\varphi(Q)$, if we let $T_P=\varphi^{-1}(T_R)$ and $T_Q=\varphi^{-1}(T_{\varphi(Q)})$, we have $T_P|_C=(d-1)P+Q$ and $T_Q|_C=dQ$. Hence, $H(\tilde{P})$ coincides with the Weierstrass semigroup as in Theorem 1.2 (b) (i).

\smallskip

(c) We consider the case where $R\neq (1:0:0)$ and $R\neq (0:0:1)$.  If we let $\varphi^{-1}(1:0:0)=\{Q\}$, $T_P=\varphi^{-1}(T_{\varphi(P)})$, and $T_Q=\varphi^{-1}(T_{\varphi(Q)})$, by the same reason as above, $H(\tilde{P})$ coincides with the Weierstrass semigroup as in Theorem 1.2 (c).

$\;$

\noindent {\bf{Acknowledgements}} 

\noindent I would like to thank Jiryo Komeda for the helpful comments concerning this work.

\end{document}